\documentclass[12pt]{article}

\usepackage[indentafter]{titlesec}
\titleformat{name=\section}{}{\thetitle.}{0.8em}{\centering\scshape}
\titleformat{name=\subsection}[runin]{}{\thetitle.}{0.5em}{\bfseries}[.]
\titleformat{name=\subsubsection}[runin]{}{\thetitle.}{0.5em}{\itshape}[.]
\titleformat{name=\paragraph,numberless}[runin]{}{}{0em}{}[.]
\titlespacing{\paragraph}{0em}{0em}{0.5em}
\titleformat{name=\subparagraph,numberless}[runin]{}{}{0em}{}[.]
\titlespacing{\subparagraph}{0em}{0em}{0.5em}

\usepackage{url}
\makeatletter
\newcommand{\address}[1]{\gdef\@address{#1}}
\newcommand{\email}[1]{\gdef\@email{\url{#1}}}
\newcommand{\sites}[1]{\gdef\@sites{\url{#1}}}
\newcommand{\@endstuff}{\par\vspace{\baselineskip}\noindent\small
	\begin{tabular}{@{}l}
		\scshape \@address \\
		\textit{E-mail address:} \@email \\
		\textit{URL:} \@sites \\[2mm]
\end{tabular}}
\AtEndDocument{\@endstuff}
\makeatother

\usepackage{titling}
\title{\normalsize\textbf{{\large O}N {\large F}INITE {\large E}XTENSIONS OF {\large L}AMPLIGHTER {\large G}ROUPS}} 
\author{Corentin Bodart}
\date{\today}

\address{Mathematical Institute, University of Oxford, UK}
\email{corentin.bodart@maths.ox.ac.uk}
\sites{https://sites.google.com/corentin-bodart}

\usepackage{import}
\usepackage{amsmath, amsfonts, amsthm, amssymb}
\usepackage{mathtools}
\usepackage{hyperref}
\usepackage[noabbrev,nameinlink,capitalise]{cleveref}
\usepackage{graphicx}
\usepackage{caption}
\usepackage{parskip}
\usepackage{adjustbox}
\usepackage{changepage}

\usepackage{tikz-cd}
\usepackage{tikz}
\usetikzlibrary{shapes, calc, arrows, patterns, decorations.markings, decorations.pathreplacing, trees}
\usepackage{xifthen}
\pgfdeclarelayer{background} 
\pgfdeclarelayer{foreground}
\pgfsetlayers{background,main,foreground}

\usepackage{easy-todo}

\usepackage{enumitem}
\usepackage{algorithm}
\usepackage{algpseudocode}

\theoremstyle{plain}
\newtheorem{thm}{Theorem}
\newtheorem*{thm*}{Theorem}
\newtheorem{lemma}[thm]{Lemma}
\newtheorem*{lemma*}{Lemma}
\newtheorem{cor}[thm]{Corollary}
\newtheorem*{cor*}{Corollary}
\newtheorem{prop}[thm]{Proposition}
\newtheorem*{prop*}{Proposition}

\theoremstyle{definition}

\newtheorem*{defi*}{Definition}
\newtheorem{rem}[thm]{Remark}
\newtheorem{conj}[thm]{Conjecture}
\newtheorem{ques}[thm]{Question}

\newcommand{\Aut}{\mathrm{Aut}}

\newcommand{\Zc}{\mathcal Z} 

\DeclareMathOperator{\supp}{supp}

\newcommand{\say}[1]{``#1"}
\newcommand{\bfe}{\mathbf e}
\newcommand{\bfu}{\mathbf u}
\newcommand{\bfv}{\mathbf v}

\newcommand{\Cay}{\mathcal C{ay}}

\newcommand{\Sc}{\mathcal S}
\newcommand{\Tc}{\mathcal T} 
\DeclareMathOperator{\Conv}{Conv}

\newcommand{\coWP}{\mathrm{coWP}}
\newcommand{\ConjGeo}{\mathrm{ConjGeo}}
\newcommand{\ot}{\leftarrow}

\newcommand{\norm}[1]{\left\|#1\right\|}
\newcommand{\abs}[1]{\left|#1\right|}
\newcommand{\la}{\left\langle}
\newcommand{\ra}{\right\rangle}
\newcommand{\lla}{\left\langle\!\left\langle}
\newcommand{\rra}{\right\rangle\!\right\rangle}
\renewcommand{\ge}{\geqslant}
\renewcommand{\le}{\leqslant}

\newcommand{\longto}{\longrightarrow}
\newcommand{\onto}{\twoheadrightarrow}
\newcommand{\into}{\hookrightarrow}
\newcommand{\acts}{\curvearrowright}

\newcommand{\N}{\mathbb N}
\newcommand{\Z}{\mathbb Z}

\newcommand{\F}{\mathbb F} 

\usepackage[backend=biber, style=numeric, sorting=nyt, 
doi=false,isbn=false,url=false,eprint=false, maxnames=50]{biblatex}
\begin{document}

\maketitle

\begin{abstract}
	We study a family of groups consisting of the simplest extensions of lamplighter groups. We use these groups to answer multiple open questions in combinatorial group theory, providing groups that exhibit various combinations of properties: \smallskip
	\begin{itemize}[leftmargin=6mm, rightmargin=2mm]
		\item Decidable Subgroup Membership and undecidable Uniform Subgroup Membership Problem.
		\item Rational volume growth series and undecidable Word Problem.
		\item Recursive (even context-free) language of conjugacy geodesics, decidable Word Problem, and undecidable Conjugacy Problem.
	\end{itemize}
	\smallskip
	We also consider the co-Word Problem, residual finiteness and the Isomorphism Problem within this class.
\end{abstract}

We introduce a family of groups. The construction starts with a group $H$, and a subset $I\subset H$ satisfying $I=I^{-1}$ and $I\not\ni 1_H$. 
\[
G(H,I) = \la a,H,z \;\middle|\; a^2=z^2=[a,z]=[h,z]=1,\;
[a,a^h]= \begin{cases} 1\!\!\! & \text{if }h\notin I \\ z\!\!\! & \text{if }h\in I \end{cases} \ra.
\]
Quotienting by $z$, we recover a presentation of the lamplighter group $C_2\wr H$. Moreover $z$ is central hence $\lla z\rra=\la z\ra\simeq C_2$ (see Proposition \ref{lem:NHI_explicit}). It follows that $G(H,I)$ is actually a central extension
\[ 1 \longto \la z\ra \longto G(H,I) \overset\tau\longto C_2\wr H \longto 1. \]
A case of special interest is $H=\Z$, in which case we denote $G(\Z,I)=G_I$. These groups were introduced by A. Genevois \cite{fournier2023no, 394830}. This construction offers a lot of flexibility through the choice of $I$. For instance, one can show that the $G_I$ form a continuum of isometric but non-isomorphic groups \cite{fournier2023no}. At the same time, they remain close to $C_2\wr\Z$ which enjoy many nice properties: it is metabelian and LERF, its Cayley graph is a Diestel--Leader graph. In turn, these properties imply decidability results for the Word, Conjugacy, and Subgroup Membership problems, and rationality of its growth series.

The construction should be compared with construction of P. Hall \cite[\S 3]{Hall}, see also \cite[\S3]{Erschler}. The groups $G_I$ are quotients of P. Hall's groups. We argue that these groups should be part of the toolbox of combinatorial group theorists trying to construct groups with pathological properties, together with $F_2\times F_2$, the Grigorchuk group, Tarski monsters and other small cancellation groups. We showcase this by answering three open questions within three different topics of geometric and combinatorial group theory.

\subsection{Subgroup Membership} We consider a group $G$, endowed with a finite generating set $S$. For convenience, we suppose that $S=S^{-1}$.
	\begin{itemize}[leftmargin=5mm]
		\item Given a subset $H\subseteq G$, we say that \emph{Membership in $H$ is decidable} if there exists an algorithm with the following specifications:
		\begin{itemize}[align=parleft, leftmargin=22mm, labelsep=12mm]
			\item[{\normalfont\textbf{Input:}}] A word $w\in S^\star$
			\item[{\normalfont\textbf{Output:}}] Yes or No, depending if $\bar w\in H$ or not.
		\end{itemize}
		\item A group $G$ has \emph{decidable Subgroup Membership} if Membership in $H$ is decidable for any finitely generated subgroup $H$.
		\item A group $G$ has \emph{decidable Uniform Subgroup Membership} if there exists an algorithm with the following specifications:
		\begin{itemize}[align=parleft, leftmargin=22mm, labelsep=12mm]
			\item[{\normalfont\textbf{Input:}}] Words $w,w_1,w_2,\ldots,w_k\in S^\star$
			\item[{\normalfont\textbf{Output:}}] Yes or No, depending if $\bar w\in \la w_1,\ldots,w_k\ra$ or not.
		\end{itemize}
	\end{itemize}
In the literature, these problems are also called \emph{Generalised Word Problem} (or the antiquated \emph{Occurrence Problem}). Most positive results solve the Uniform Subgroup Membership. This includes the results for metabelian groups, and finitely presented LERF groups.  In the other direction, most undecidability results only look at a fixed subgroup. Most notably, Mikhailova proves that membership in the fixed subgroup
\[ M(G) = \{(u,v)\in F_2 \times F_2 \mid \psi(u)= \psi(v) \} \le F_2\times F_2 \]
is undecidable \cite{Mihailova}, where $\psi\colon F_2\onto G$ is a morphism onto a finitely presented group with undecidable Word Problem. This explains why it is hard to separate both problems. We provide an example within the family $G_I$:

\begin{thm}
	There exists a recursive subset $I\subset\Z$ such that $G_I$ has decidable Subgroup Membership and undecidable Uniform Subgroup Membership.
\end{thm}

Boiling the proof down to a few lines, finitely generated subgroups $H$ of $G_I$ fall into two classes, with an algorithm deciding membership uniformly in either class. However we cannot compute in which class a subgroup falls given a generating set, hence the Uniform Subgroup Membership is undecidable.

\subsection{Growth series and the Word Problem} Starting from the same data $(G,S)$ as in the previous paragraph, we turn $G$ into a metric space and define a few geometric notions:
\begin{itemize}[leftmargin=8mm]
	\item The \emph{word length} of an element $g\in G$ is
	\[\norm{g}_S\coloneqq \min\bigl\{n\ge 0 \mid \exists w\in S^n,\; \bar w=g\bigr\}. \]
	\item The \emph{volume growth} of $(G,S)$ is the function
	\[ \beta_{G,S}(n) \coloneqq \#\bigl\{g\in G\mid \norm{g}_S\le n\bigr\}. \]
	\item The \emph{volume growth series} of $(G,S)$ is the formal series
	\[ \Gamma_{G,S}(x) \coloneqq \sum_{n\ge 0} \beta_{G,S}(n)\cdot x^n \in \Z[[x]]. \]
\end{itemize}

Historically, computing growth series is strongly linked to languages in groups, eg.\ normal forms and the Word Problem. For instance, computations for some Fuchsian groups in \cite{cannon1984combinatorial} led to the definition of automatic groups. Some landmark results state that the growth series of virtually abelian groups and hyperbolic groups is always rational. In the opposite direction, we know of very few negative results. One of these is an observation of Thurston:
\begin{thm*}[{Thurston, see \cite{cannon1980growth, grigorchuk1997problems}}]
	Suppose that $G$ is a recursively presented group with undecidable Word Problem, then the growth series $\Gamma_{G,S}(x)$ is not computable (in particular not rational nor algebraic).
\end{thm*}
We prove that the hypothesis \say{recursively presented} cannot be removed, answering a question of Duchin and Shapiro \cite[p.221]{Heis_is_panrational}.
\begin{thm}
	There exists a pair $(G_I,S)$ such that $G_I$ has undecidable Word Problem and the growth series $\Gamma_{G_I,S}(x)$ is rational.
\end{thm}
More strongly, we use that $G_I$ is isometric to $G_\emptyset=C_2\times (C_2\wr\Z)$ (for the correct choice of generating sets). This fact was observed by Genevois \cite{394830}.

\subsection{Conjugacy geodesics and the Conjugacy Problem} Finally, we recall a few constructions related to conjugacy:
\begin{itemize}[leftmargin=8mm]
	\item Two elements $g,h\in G$ are \emph{conjugated} if there exists $c\in G$ such that $g=cgc^{-1}$. In that case, we write $g\sim h$ (or $g\sim_Gh$ when ambiguous).
	\item Given $g\in G$, we define $\norm{[g]}_S\coloneqq \min\bigl\{\norm h_S \;\big|\; h\sim g\bigr\}$.
	\item The language of \emph{conjugacy geodesics} of $(G,S)$ is
	\[ \ConjGeo(G,S) = \bigl\{w\in S^\star \;\big|\; \ell(w)=\norm{[\bar w]}_S\bigr\}.\]
\end{itemize}
In many cases, strong geometric properties leading to algorithms solving the Conjugacy Problem also implies that $\ConjGeo(G,S)$ is rational \cite{ConjGeo_intro,ConjGeo_main, ConjGeo_weakhyp, ConjGeo_virtRAAG, ConjGeo_dihedral}. More generally, if the Conjugacy Problem is decidable, then the language $\ConjGeo(G,S)$ is recursive. The other way around, the fact that $\ConjGeo(G,S)$ is low in Chomsky's hierarchy does not easily imply that the Conjugacy Problem is decidable. The possibility that this implication might fail was raised by Ciobanu, Hermiller, Holt and Rees \cite{ConjGeo_main}. They even suggest a pair $(G,S)$ with undecidable Conjugacy Problem and recursive $\ConjGeo(G,S)$. We confirm their suspicion in a stronger sense:
\begin{thm}
	There exists $(G,S)$ such that $\ConjGeo(G,S)$ is context-free, the Word Problem is decidable, and the Conjugacy Problem is undecidable.
\end{thm}
More precisely, we take $G=G(F_2,I)$ for some recursive subset $I\subset F_2$. The relatively low complexity of $\ConjGeo(G,S)$ comes from the quotient $C_2\wr F_2$. 

\subsection*{Acknowledgment} I would like to thank Carl-Fredrik Nyberg-Brodda for bringing the issue of uniform vs non-uniform subgroup membership. The author is supported by a Postdoc.Mobility grant from the Swiss SNF.
\counterwithin{thm}{section}
\section{Generalities}

Using Tietze moves, we can see that $G(H,I)=N(H,I)\rtimes H$, where
\[ N(H,I) = \la \{a_h\}_{h\in H}, z \;\middle|\; a_h^2=z^2=[a_h,z]=1,\; [a_g,a_h]=\begin{cases} 1 \!\! & \text{if }g^{-1}h\notin I \\ z \!\! & \text{if }g^{-1}h\in I \end{cases}\ra. \]
Let us \say{realise} $N(H,I)$ as an explicit $2$-step nilpotent group.
\begin{prop} \label{lem:NHI_explicit}
	Fix a total order $<$ on $H$. Then $N(H,I)$ is isomorphic to $(\bigoplus_H\F_2)\times \F_2$ with product $(\bfu,m)(\bfv,n) = \big(\bfu+\bfv,\; m+n+\omega_I(\bfu,\bfv)\bigr)$, where
	\[ \omega_I(\bfu,\bfv) = \sum_{g<h} \chi_I\bigl(g^{-1}h\bigr)\,\bfu_g\bfv_h \pmod2.\]
	In particular, we see that $z=(\mathbf 0,1)\ne e$. Moreover $N(H,I)$ is $2$-step nilpotent of exponent $4$, hence $N(H,I)$ is locally finite.
\end{prop}
\begin{proof}
	One can check that the product is associative using that $\omega_I$ is bilinear, and $(\mathbf 0,0)$ is a neutral element. Moreover every element is invertible since
	\[ (\bfu,m)^4=(\mathbf 0,\omega_I(\bfu,\bfu))^2=(\mathbf 0,0), \]
	so $N\coloneqq (\bigoplus_H\F_2)\times \F_2$ is a group! The morphism $F(a_h,z)\to N$ defined by $a_h\mapsto (\delta_h,0)$ and $z\mapsto (\mathbf 0,1)$ sends relations to $(\mathbf 0,0)$ using that $I=I^{-1}$ and $I\not\ni 1$, so the morphism descends to a map $N(H,I)\onto N$. Finally, every element of $N(H,I)$ can be written as a product $$a_{h_1}\ldots a_{h_\ell}z^\varepsilon$$ with $h_1>\ldots >h_\ell$ and $\varepsilon\in\{0,1\}$ using the relations, and these are mapped to distinct elements of $N$ (specifically, $(\delta_{h_1}+\ldots+\delta_{h_\ell}, \varepsilon)$), so the epimorphism $N(H,I)\onto N$ is injective: it is an isomorphism.
\end{proof}

\medskip

The next two results follow easily:
\begin{cor} \label{cor:commutator}
	Given an element $x=(\bfu,m)\in N(H,I)$, we define its support as $\supp(x)=\{g\in H:\bfu_g\ne 0\}$. Then $[x,y] = z^\varepsilon$ where
	\[ \varepsilon=\sum_{g\in\supp(x)}\sum_{h\in\supp(y)} \chi_I(g^{-1}h) \pmod2.\]
\end{cor}

\begin{cor} \label{lem:natural} Consider two groups $H,H'$ and $I\subset H'$.
	\begin{enumerate}[leftmargin=8mm, label={\normalfont(\alph*)}]
		\item Any morphism  $\iota\colon H'\into H$ induces a morphism $G(H',I)\into G(H,\iota(I))$.
		\item Any morphism $\psi\colon H\onto H'$ induces a morphism $G(H,\psi^{-1}(I))\onto G(H',I)$.
	\end{enumerate}
\end{cor}

Finally, we can characterise when $G(H,I)$ has decidable Word Problem:
\begin{thm} \label{thm:WP}
	$G(H,I)$ has decidable Word Problem if and only if $H$ has decidable Word Problem and $I$ is recursive.
\end{thm}
\begin{proof}
	We first prove the \say{if} direction. Consider a word $w\in (S_H\cup\{a,z\})^*$. Since $H$ has decidable Word Problem, one can decide if $\tau(\bar w)=1$ or not. In the latter case, one conclude that $\bar w\ne 1$. In the former, we can write
	\[ \bar w = a_{\pi(w_1)} a_{\pi(w_2)} \ldots a_{\pi(w_n)} z^\varepsilon \] 
	where $w_i$ is prefix of $w$ before the $i$-th use of $a$, and $\varepsilon\in\{0,1\}$ counts the parity of the number of $z$ in $w$. Moreover, we can form pairs of $a_{\pi(w_i)}$ that cancel out in $C_2\wr H$. Now we just have move the factors around, while keeping track of the extra $z$ factors if one need to permute $a_{\pi(w_i)}$ and $a_{\pi(w_j)}$ where $\pi(w_i)^{-1}\pi(w_j)\in I$ (which is decidable since $I$ is recursive).
	
	For the other direction, note that $H$ is a subgroup of $G(H,I)=N(H,I)\rtimes H$ hence $H$ has decidable Word Problem. Moreover, we can decide whether $[a,hah^{-1}]z=1$ or not, hence whether $h\in I$ or not.
\end{proof}

\begin{center}
	\adjustbox{scale=.95}{%
		\begin{tikzcd}
			& & 1\ar[d] & 1\ar[d] & \\
			1 \ar[r] & \la z\ra \ar[d,"\wr"]\ar[r, hook] & N(H,I) \ar[r, two heads, "\tau"]\ar[d, hook] & \bigoplus_H \F_2 \ar[d, hook]\ar[r] & 1 \\
			1 \ar[r] & \la z\ra \ar[r, hook] & G(H,I) \ar[r, two heads, "\tau"]\ar[d, two heads, swap, "\pi"] & \F_2\wr H \ar[r]\ar[d, two heads] & 1 \\
			&  & H \ar[r, "\sim"] \ar[d]\ar[u, bend right] &  H \ar[d] \ar[u, bend right] & \\
			& & 1 & 1 &
	\end{tikzcd}}
	\captionsetup{font=small}
	\captionof{figure}{The third isomorphism theorem}
\end{center}

\section{Membership problems}

We construct a finitely generated group with decidable Subgroup Membership, and undecidable Uniform Subgroup Membership. More strongly,
\begin{thm} \:
	\begin{enumerate}[leftmargin=8mm, label={\normalfont(\alph*)}]
		\item $G_I$ has decidable Subgroup Membership if $I\subset\Z$ is recursive.
		\item There exists a recursive set $I\subset \Z$ such that the decision problem
		\begin{itemize}[align=parleft, leftmargin=22mm, labelsep=12mm]
			\item[{\normalfont\textbf{Input:}}] Two elements $g_1,g_2\in G_I$ {\normalfont(}given as words over $\{a,t^\pm,z\})$.
			\item[{\normalfont\textbf{Output:}}] Yes or No, depending if $z\in\la g_1,g_2\ra$ or not.
		\end{itemize}
		is undecidable.
	\end{enumerate}
\end{thm}
This last problem is an instance of \say{fixed-target} membership problem, as introduced in \cite{Gray_NB}. We start by proving (a) using the following result:
\begin{thm}[Romanovskii \cite{subgroup_meta}] \label{thm:subgroup_metabelian}
	If $G$ is finitely generated metabelian group, then $G$ has decidable Uniform Subgroup Membership.
\end{thm}
We fix a finitely generated subgroup $H\le G_I$ and separate two cases:
\begin{algorithm}[H]
	\caption{(Is $\bar w\in H$? Case when $z\not\in H$)}\label{alg:z_out}
	\begin{algorithmic}[1]
		\Require A word $w\in \{a,t^\pm,z\}^\star$
		\If{$\pi(\bar w)\in\pi(H)$} \Comment{Theorem \ref{thm:subgroup_metabelian}}
			\State Find $h\in H$ such that $\pi(\bar w)=\pi(h)$. 
			\If{$\bar w=h$} \Comment{Theorem \ref{thm:WP}}
				\State\Return\say{Yes, $\bar w\in H$}
			\Else
				\State\Return\say{No, $\bar w\notin H$}
			\EndIf
		\Else
			\State\Return\say{No, $\bar w\notin H$}
		\EndIf
	\end{algorithmic}
\end{algorithm}
On line 2, we use that $H$ is finitely generated to enumerate elements $h\in H$, until we eventually find an element such that $\pi(\bar w)=\pi(h)$. The answer on line 6 comes from the fact that $h\in H$ and $z\notin H$ implies $hz\notin H$.

\smallskip

When $z\in H$, we can use the following even simpler algorithm:
\begin{algorithm}[H]
	\caption{(Is $\bar w\in H$? Case when $z\in H$)}\label{alg:z_in}
	\begin{algorithmic}[1]
		\Require A word $w\in \{a,t^\pm,z\}^\star$
		\If{$\pi(\bar w)\in\pi(H)$} \Comment{Theorem \ref{thm:subgroup_metabelian}}
			\State\Return\say{Yes, $\bar w\in H$}
		\Else
			\State\Return\say{No, $\bar w\notin H$}
		\EndIf
	\end{algorithmic}
\end{algorithm}

Finally, we find a recursive subset $I\subset\Z$ such that $G_I$ has undecidable Uniform Subgroup Membership. We start with the following lemma:
\begin{lemma} \label{lem:Turing_membership}
	There exists a recursive set $I\subset \mathbb Z$ (with $0\notin I$ and $I=-I$) such that, given $n\in\N$, one cannot decide if there exists $i\in I$ such that $n\mid i$.
\end{lemma}
\begin{proof}
	Fix a computable enumeration $\N\leftrightarrow\{\text{Turing machines}\}$, and $p_1,p_2,\ldots$ the prime numbers in increasing order. We define a set $I$ whose elements are powers of primes: $\pm p_n^m\in I$ if the $n$-th Turing machine halts after $m$ steps. We just have to observe that there exists $i\in I$ such that $p_n\mid i$ if and only if the $n$-th Turing machine halts, which is undecidable.
\end{proof}

We observe that $\la a,t^n,z\ra\simeq G_J$, where $J=\{j\in\Z\mid nj\in I\}$. If $J$ is empty, then $G_J=\langle a,t^n\rangle\times\langle z\rangle$ and $z\notin\langle a,t^n\rangle$. Otherwise $z=[a,t^{jn}at^{-jn}]\in\la a,t^n\ra$ for some $j\in J$. This allows the following reduction:
\[ z\in\la a,t^n\ra \iff J\ne\emptyset \iff \exists i\in I \text{ such that } n\mid i, \]
the last problem being undecidable by Lemma \ref*{lem:Turing_membership}. \hfill$\square$

\medskip

\begin{rem}Given a finitely generated group $G=\la S\ra$, and a finitely generated subgroup $H=\la T\ra$, we define the (upper) \emph{distortion function} as
	\[ \Delta_G^H(n) = \max\bigl\{ \norm{g}_T \;\big|\; g\in H,\, \norm{g}_S\le n \bigr\}.\]
	When $G$ has decidable Word Problem, membership in $H$ is decidable if and only if $\Delta_G^H(n)$ admits a recursive upper bound. For $G=C_2\wr\Z$, Davis and Ol'shanskii proved that $\Delta_G^H(n)$ is always linear \cite{distortion_ZwrZ}. This extends to $G_I$:
	
	We fix generating sets $\{a,az,t^\pm,t^\pm z\}$ and $\{a,t^\pm\}$ for $G_I$ and $C_2\wr\Z$, so that $\norm{g}=\norm{\tau(g)}$ for all $g\ne z$. Consider $H=\la T\ra\le G_I$. There are two cases:
	\begin{itemize}[leftmargin=8mm]
		\item If $z\notin H$, then $\norm{h}_T=\norm{\tau(h)}_{\tau(T)}$ for all $h\in H$, hence
		\[ \Delta_{G_I}^H(n) = \Delta_{C_2\wr \Z}^{\tau(H)}(n) \preceq n.\]
		\item If $z\in H$, then $\norm{h}_T\le\norm{\tau(h)}_{\tau(T)}+\norm{z}_T$ for all $h\in H$, hence
		\[ \Delta_{G_I}^{H}(n) \le \Delta_{C_2\wr\Z}^{\tau(H)}(n) + \norm{z}_T \preceq n.\]
	\end{itemize}
	This doesn't contradict the fact that the Uniform Subgroup Membership is undecidable for some $I$. Indeed, in this case, the implicit constants in the bound $\Delta_{G_I}^{H}(n)\preceq n$ cannot be computed effectively from $T$.
\end{rem}

\subsection{Other membership problems} We observe that $G_I$ has decidable Submonoid Membership as soon as $I$ is recursive:
\begin{prop}
	Let $(G,S)$ be a group given as an extension
	\[ 1 \longto N \longto G \overset\pi\longto \Z \longto 1 \]
	with $N$ locally finite. Then (Uniform) Subgroup Membership and (Uniform) Semigroup Membership are Turing equivalent.
\end{prop}
\begin{proof}
	Fix a f.g.\ submonoid $M=\{g_1,g_2,\ldots,g_n\}^*$. There are two cases:
	\begin{itemize}[leftmargin=8mm]
		\item If there exist $i,j$ with $\pi(g_i)<0<\pi(g_j)$, then $M$ is a subgroup of $G$. Indeed, for all $g\in M$, there exists $p,q,r\in\N$ such that
		\[ \pi(g^pg_i^qg_j^r) = p\pi(g)+q\pi(g_i)+r\pi(g_j)=0,\] 
		and $\ker\pi=N$ is locally finite, hence $g$ is invertible in $M$.
		\item Otherwise, given $w\in S$, any product of elements $g_i$ evaluating to $\bar w$ may only contain at most $\abs{\pi(\bar w)}$ factor $g_i$ with $\pi(g_i)\ne 0$. Moreover $H\coloneqq\{g_i\mid \pi(g_i)=0\}^*\subseteq N$ is a finite subgroup which can be enumerated given an oracle for the Word Problem. It follows that there is only a finite list of products which could potentially equal $\bar w$:
		\[ \left\{h_0k_1h_1k_2\ldots k_\ell h_\ell \;\Big|\; h_j\in H,\, k_j\in\{g_i\mid \pi(g_i)\ne 0\}\text{ and }\ell\le\abs{\pi(\bar w)}\right\}. \]
		$\bar w\in M$ if and only if one of these product equals $\bar w$.
	\end{itemize}
	To summarize, the are two cases and we can effectively determine in which case we fall. In the first case, $M=\la g_1,\ldots,g_n\ra$ is a subgroup. In the second case, Membership in $M$ is Turing reducible to the Word Problem.
\end{proof}
It follows that the same group $G_I$ also separates Submonoid Membership and Uniform Submonoid Membership. This leaves open the following question:
\begin{ques}
	Does $G_I$ have decidable Rational Subset Membership as soon as $I$ is recursive? Note that the problem is decidable in $C_2\wr\Z$ \cite{rat_wreath}.
\end{ques}
If the answer is negative, this would separate (non-uniform) Semigroup and Rational Subset Membership, giving a second example after \cite{Dong_submonoid}.

Another intriguing possibility would be to separate Uniform Submonoid Membership and the Knapsack problem. For $I$ selected at random, it seems that any finitely generated group $H\le G_I$ should fall in one of two cases:
\begin{itemize}[leftmargin=8mm]
	\item Either $H\le N(\Z,I)$ is finite,
	\item Or $H$ is infinite and $H\ni z$.
\end{itemize}
Therefore Uniform Subgroup Membership should be decidable for $I$ recursive but sufficiently irregular. This opens the following possibility:
\begin{ques}
	Does there exists $I\subset\Z$ such that
	\begin{enumerate}[leftmargin=8mm, label=(\arabic*)]
		\item the Uniform Subgroup Membership problem is decidable in $G_I$, and
		\item the (Uniform) Knapsack problem is undecidable in $G_I$?
	\end{enumerate}
\end{ques}

\section{Growth series and the Word Problem}

We first recall and prove an observation due to Thurston.
\begin{prop}[{Thurston, see \cite{cannon1980growth, grigorchuk1997problems}}] \label{prop:Thurston}
    Suppose that $(G,S)$ is recursively presented and $G$ has undecidable Word Problem, then the growth series $\Gamma_{G,S}(x)$ is not computable (in particular not rational nor algebraic).
\end{prop}
\begin{proof}
	We prove the contrapositive: we suppose that the growth series is computable and provide an algorithm solving the Word Problem. Since $(G,S)$ is recursively presented, there exists a sub-routine enumerating all relations $u_t=v_t$ between words $u_t,v_t\in S^\star$.
	\begin{algorithm}[H]
		\caption{(The Word Problem in $(G,S)$ from its growth series)}
		\begin{algorithmic}[1]
			\Require Two words $u,v\in S^\star$
			\State $n\gets \max\{\abs{u},\abs{v}\}$
			\State $\beta\gets \beta_{G,S}(n)$ \Comment{$\Gamma_{G,S}(x)$ is computable}
			\medskip
			\State $t\gets 0$
			\State $\sim \;\gets$ equality relation on $S^{\le n}$.
			\While{ $\abs{S^{\le n}/\!\sim\,}>\beta$ }
			\State $t\gets t+1$
			\State $\sim \;\gets$ transitive closure of $\sim$ and $u_t=v_t$ (if $\abs{u_t},\abs{v_t}\le n$).
			\EndWhile
			\If{$u\sim v$}
			\State\Return{\say{Yes, $u=v$ in $G$}}
			\Else
			\State\Return{\say{No, $u\ne v$ in $G$}}
			\EndIf
		\end{algorithmic}
	\end{algorithm} \vspace*{-3mm}
	For any recursively presented group, the algorithm enumerating relations will eventually enumerate all relations between words of length at most $n$. Knowing $\beta_{G,S}(n)$ gives a certificate when all such relations have been given, allowing the algorithm to terminate.
\end{proof}

\subsection{Example}
We prove that the hypothesis \say{recursively presented} in Proposition \ref{prop:Thurston} cannot be removed, answering a question of Duchin and Shapiro \cite[p.221]{Heis_is_panrational}, and correcting a recurring mistake \cite[p.191]{Putnam}, \cite[p.181]{Choi_Ho_Pengitore_2022}.

\begin{prop}
Consider the group $G_I$ together with the generating set $S=\{a,az,t^\pm,t^\pm z\}$. Then the growth series $\Gamma_{G_I,S}(x)$ is rational.
\end{prop}
\begin{proof} Note that $S=\tau^{-1}(T)$ where $T=\{a,t^\pm\}$ is the standard generating set of $C_2\wr \Z$. It follows that $\norm{g}_S = \norm{\tau(g)}_T$ as soon as $g\ne \{1,z\}$. Combined with $\norm{1}_S=0$ and $\norm{z}_S=2$, this implies that
\[ \Gamma_{G_I,S}(x) = 1+x^2+ \abs{\ker\tau}\cdot  \big( \Gamma_{C_2\wr\Z,T}(x)-1 \big) \]
or more explicitly
\[ \Gamma_{G_I,S}(x) = 1 + x^2 + 2 \left( \frac{(1+x)^3(1-x)^2(1+x+x^2)}{(1-x^2-x^3)^2(1-x-x^2)}-1 \right) \]
using computations from \cite{Johnson}.
\end{proof}

Note that this is does not extends to any generating set:
\begin{lemma} \label{lem:marked_groups}
	Consider $S'=\{a,t^\pm,z\}$ and $I,J\subset \Z$ such that $I\cap [1,r]=J\cap[1,r]$, then balls of radius $2r+3$ in the labelled Cayley graphs $\Cay(G_I,S')$ and $\Cay(G_J,S')$ are isomorphic. In particular, the map
	\[ \begin{pmatrix}
		\bigl\{ I\subset\Z \;\big|\; I=-I,\; I\not\ni 0\bigr\} & \longto & \mathcal M \\
		I & \longmapsto & (G_I,S')
	\end{pmatrix}\]
	is continuous, where $\mathcal M$ is the space of marked groups.
\end{lemma}
\begin{proof}Suppose that we have two words $u,v$ over $S'$ such that $\bar u=\bar v$ in $G_I$ and $\bar u\ne \bar v$ in $G_J$. In particular, we have $\tau(\bar u)=\tau(\bar v)$, and there exist two positions $m,n\in\Z$ such that $m-n\in I\Delta J$ and
\begin{itemize}[leftmargin=8mm]
	\item $u$ switch on the lamps at position $m$ and $n$ in that order, and
	\item $v$ switch on the lamps at position $n$ and $m$ in that order.
\end{itemize}
(Possibly, the lamps get turned on and off at other points along $u$ and $v$.) Consider $p=\pi(\bar u)=\pi(\bar v)$ the end position of both words. This forces that
	\begin{align*}
		\ell(u)-2 & \ge d(0,m) + d(m,n) + d(n,p), \\
		\ell(v)-2 & \ge d(0,n) + d(n,m) + d(m,p).
	\end{align*}
	Depending on the relative positions of $0$, $m$, $n$ and $p$ on the line, at least one of the two words has length greater than $2d(m,n)+2\ge 2r+4$.
\end{proof}
\begin{prop} Consider the generating set $S'=\{a,t^\pm,z\}$. Then the growth series $\Gamma_{G_I,S'}(x)$ is computable if and only if $I$ is recursive.
\end{prop}
\begin{proof}
    Suppose that $\Gamma_{G_I,S'}$ is computable, i.e., there exists an algorithm which outputs the values $\beta_{G_I,S'}(n)$. We compute $I_r=I\cap [-r,r]$ recursively. First $I_0=\emptyset$. Suppose we have computed $I_r$ for some $r\ge 0$, then
    \[
    \beta_{G_I,S'}(2r+4) = \begin{cases}
    	\beta_{G_{I_r},S'}(2r+4) & \text{if }r+1\notin I,\\
    	\beta_{G_{I_r},S'}(2r+4)-r-2 & \text{if }r+1\in I.
    \end{cases}
    \]
    Indeed, most of the two balls are isomorphic by Lemma \ref{lem:marked_groups}. The only words that represent distinct elements in $G_{I_r}$ and the same element in $G_{I_r\cup\{r+1\}}$ are $t^{-j}at^iat^{j-i}$ and $t^{i-j}at^{-i}at^j$ for each $j\in\{0,\ldots,r+1\}$. Since $\beta_{G_{I_r},S'}(2r+4)$ is computable by induction, we can compute $I_{r+1}$.
\end{proof}
This improves on a result of Stoll proving the existence of a group such that the growth series is rational w.r.t.\ a generating set, and transcendental w.r.t.\ another generating set \cite{Stoll}. So far, his computations for $H_5(\Z)$ with \say{standard} and \say{cubical} generating sets was the only source of examples. (Analogous results can be derived for many groups from this first example, for instance higher Heisenberg groups $H_n(\Z)$ for $n\ge 5$ odd, or $H_5(\Z)*\Z$.)
\begin{ques}
	Find necessary and sufficient conditions on $I$ ensuring that $\Gamma_{G_I,S'}(x)$ is rational. Is the series rational as soon as $I$ is eventually periodic?
\end{ques}
\section{Conjugacy Geodesics and the Conjugacy Problem}
In this section, we prove the following theorem:
\begin{thm} \label{thm:bad_conjugacy}
	There exists a group $G$ generated by a finite set $S$ such that
	\begin{itemize}[leftmargin=8mm]
		\item $G$ has decidable Word Problem,
		\item $G$ has undecidable Conjugacy Problem, and
		\item $\ConjGeo(G,S)$ is unambiguously context-free.
	\end{itemize}
\end{thm}
More specifically, we take $G=G(F_2,I)$ for some recursive subset $I\subset F_2$. In particular, this ensures that the Word Problem is decidable (Theorem \ref{thm:WP}). The other two points are proven over the next two subsections.

\subsection{Undecidability} We start the proof by a lemma:

\begin{lemma} \label{lem:conj_to_sum}
	Let $H$ be a torsion-free group and $x\in N(H,I)$. Then
	\[ x\sim xz \text{ in }G(H,I) \iff \exists g\in H,\; \abs{g^{-1}\cdot \supp(x)\cap I} \text{ is odd.}\]
\end{lemma}
\begin{proof}
	First, if $x\in\{1,z\}$ is central, then $x\not\sim xz$ and $\supp(x)=\emptyset$, so the statement trivially holds. We suppose that $x\notin\{1,z\}$, i.e., $\supp(x)\ne\emptyset$.
	
	Suppose that $x\sim xz$ and let $c\in G(H,I)$ be a conjugator. We have
	\[ cxc^{-1}=hz \implies \pi(c)\cdot \supp(x)=\supp(cxc^{-1})=\supp(x) \]
	hence $\pi(c)=1_H$ since $\supp(x)\ne\emptyset$ and $H$ is torsion-free. Since $c\in N(H,I)$, we can use Corollary \ref{cor:commutator} to compute $[c,x]=z^\varepsilon$ with
	\[ \varepsilon = \sum_{g\in \supp(c)}\sum_{h\in\supp(x)} \chi_I(g^{-1}h) = \sum_{g\in\supp(c)} \abs{g^{-1}\cdot \supp(x)\cap I} \pmod 2.\]
	Since $cxc^{-1}=xz$, we have $\varepsilon\equiv 1$ hence there exists $g\in H$ such that the intersection has odd size. Reciprocally, if such a $g\in H$ exists, we can take $c=gag^{-1}\in N(H,I)$ which satisfies $\supp(c)=\{g\}$ hence $cxc^{-1}=xz$.
\end{proof}
This motivates the next lemma:
\begin{lemma} \label{lem:Turing_conjugate}
	There exists a recursive subset $J\subset\bigoplus_\Z\Z$ (satisfying $\mathbf 0\notin J$ and $J=-J$) such that, given a finite subset $X\subset\bigoplus_\Z\Z$, it is not decidable whether there exists $g\in\bigoplus_\Z\Z$ such that $\abs{(X-g)\cap J}$ is odd.
\end{lemma}
\begin{proof}
	Fix a computable enumeration $\Z\leftrightarrow\{\text{Turing machines}\}$, and define
	\[ J = \left\{ (x_i)_i\in\bigoplus_\Z\Z \;\middle|\;  \begin{array}{c}\exists i\in\Z \text{ such that the $i$-th Turing} \\ \text{machine halts in $\abs{x_i}$ steps or less}\end{array} \right\}. \]
	Let $\bfe_i$ be the vector with $1$ in the $i$-th coordinate and $0$ elsewhere, and $X=\{\mathbf 0,\bfe_i\}$. There exists $g$ such that $\abs{(X-g)\cap J}$ is odd if and only if $J$ is not invariant under translation by $\bfe_i$, i.e., if the $i$-th machine halts.
\end{proof}

Finally, we construct our example. We consider the epimorphism 
\[ \psi\colon F_2 = \la s,t\ra \longto \Z\wr\Z=\la s,t \;\big|\; [s,t^nst^{-n}]=1 \;(n\ge 1)\ra. \]
and define $G=G(F_2,I)$ where $I=\psi^{-1}(J)$. To each finite subset $X=\{h_1,\ldots,h_n\}\subset \bigoplus_\Z\Z$, we associate the computable element
\[ x = w_{h_1}aw_{h_1}^{-1} \cdot \ldots \cdot w_{h_n}aw_{h_n}^{-1} \in G(F_2,I),\]
where $w_h\in F_2$ satisfies $\psi(w_h)=h$. (Take for instance $w_h$ the ShortLex geodesic for $h$, seen as an element in $F_2$.) By construction $\psi$ restricts to a 1-to-1 map $\supp(x)\to X$. Using Lemma \ref{lem:conj_to_sum}, we have $x\sim xz$ if and only if 
\[ \exists g\in F_2,\; \abs{g^{-1}\cdot \supp(x)\cap I}=\abs{\psi(g)^{-1}\cdot X\cap J} \equiv 1\pmod2 \]
If $\psi(g)\notin\bigoplus_\Z\Z$, then $\psi(g)^{-1}\cdot X\cap J=\emptyset$ (since $J,X\subset\bigoplus_\Z\Z$). Therefore, the problem reduces to the existence of $g\in \psi^{-1}(\bigoplus_\Z\Z)$ such that
\[ \abs{(X-\psi(g))\cap J}\text{ is odd,} \]
which is undecidable by Lemma \ref{lem:Turing_conjugate}.
\begin{rem}
	This example proves that hypothesis (iii) in \cite[Theorem 3.1]{Martino} cannot be removed. Here, the centraliser of an element $x\in \bigoplus_{F_2}C_2$ is $C_{C_2\wr F_2}(x)=\bigoplus_{F_2}C_2$, which is much larger than $\la x\ra$.
\end{rem}

\subsection{Language of conjugacy geodesics} We prove that, contrary to decidability of the Conjugacy Problem, the complexity of $\ConjGeo(G,S)$ is preserved  by extension with finite kernel (for specific generating sets):

\begin{lemma}
	Consider a short exact sequence
	\[ 1 \longto F \longto G \overset{\tau}\longto Q \longto 1,\]
	with $F$ finite. Let $T$ be a finite generating set of $Q$, and define the associated generating set $S=\tau^{-1}(T)\cup F$ of $G$. Then
	\[ \ConjGeo(G,S) = \tau^{-1}\bigl(\ConjGeo(Q,T)\bigr) \cup (F\setminus\{1\}).\]
	In particular, given any class of languages $\mathbf C$ closed under inverse image by non-erasing homomorphisms and union with finite languages, $\ConjGeo(G,S)$ belongs to $\mathbf C$ if and only if $\ConjGeo(Q,T)$ belongs to $\mathbf C$. This covers the class of regular (resp.\ one-counter, context-free, recursive, etc.) languages.
\end{lemma}
\begin{proof}
	Consider $g\in G\setminus F$. We have $[g]\subset G\setminus F$ and $\tau([g])=[\tau(g)]$. Moreover, $\norm{h}_S=\norm{\tau(h)}_T$ for all $h\in G\setminus F$. We conclude that
	\[ \norm{[g]}_S = \min\bigl\{\norm h_S \mid h\sim g\bigr\} = \min\{\norm{k}_T \mid k\sim\tau(g)\} = \norm{[\tau(g)]}_T. \]
	 It follows that $w\in S^*$ is a conjugacy geodesic if and only if $\bar w\ne 1_Q$ and $\tau(w)\in T^*$ is a conjugacy geodesic, or $\bar w\in F$ and $w\in \{\varepsilon\}\cup (F\setminus\{ 1\})$.
\end{proof}

The fact that the language $\ConjGeo(G(F_2,I),S)$ is unambiguously context-free (for $S=\pi^{-1}(\{a,s^\pm,t^\pm\})$) now follows from Theorem \ref{thm:conjgeo_wreath}.
\begin{rem}
	The same argument shows that $G(\Z\wr\Z,I)$ has decidable Word Problem and undecidable Conjugacy Problem. Moreover, Jane Matthews has proven that $C_2\wr(\Z\wr\Z)$ (among many wreath products) has decidable Conjugacy Problem \cite{Matthews}, hence $\ConjGeo(G,S)$ is at least recursive in this case. However, the travelling salesman problem on $\Z\wr\Z$ is much harder than that on $F_2$, hence we cannot have much hope to show that $\ConjGeo(G,S)$ is any better than context-sensitive, and most likely not context-free.
\end{rem} 

\subsection{Conjugacy in $G(\Z,I)$} One may wonder why we did not use $G_I$ in the previous construction. The reason is the following result:

\begin{thm}
	If $I$ is recursive, then $G_I$ has decidable Conjugacy Problem.
\end{thm}
\begin{proof} We first describe a sub-routine to decide whether $h\sim hz$ or not. We separate two cases, depending if $I\subset\Z$ is periodic or not.
	
If $I$ is periodic, we can use the following algorithm:
\begin{algorithm}[H]
	\caption{(Is $h\sim hz$? Case when $I$ is $p$-periodic.)}\label{alg:period}
	\begin{algorithmic}[1]
		\Require{$h\in G_I$}
		\State Answer $\gets$ No
		\If{$\pi(h)=0$ \textbf{and} $\tau(h)\ne 1$} \Comment{i.e., $h\in N(\Z,I)\setminus\{1,z\}$}
			\For{$k\in \{0,\ldots,p-1\}$}
				\State $r\gets \sum_{n\in\supp(\tau(h))} \chi_I(n-k)$
				\If{$r\equiv 1\pmod 2$}
					\State Answer $\gets$ Yes
				\EndIf
			\EndFor
		\EndIf
		\State\Return Answer
	\end{algorithmic}
\end{algorithm}

If $I$ is not periodic, we can use the following even simpler algorithm:
\begin{algorithm}[H]
	\caption{(Is $h\sim hz$? Case when $I$ is not periodic.)}\label{alg:aperiod}
	\begin{algorithmic}[1]
		\Require{$h\in G_I$}
		\State Answer $\gets$ No
		\If{$\pi(h)=0$ \textbf{and} $\tau(h)\ne 1$} \Comment{i.e., $h\in N(\Z,I)\setminus\{1,z\}$}
			\State Answer $\gets$ Yes
		\EndIf
		\State\Return Answer
	\end{algorithmic}
\end{algorithm}

Both algorithms are clearly terminating, we prove correctness.
\begin{itemize}[leftmargin=8mm]
	\item If $\tau(h)=1$, i.e., $h\in\{1,z\}$, then $h$ is central hence $h\not\sim hz$.
	\item If $\tau(h)\in\bigoplus_\Z C_2\setminus\{1\}$. By Lemma \ref{lem:conj_to_sum}, we have
	\[ h\sim hz \iff \exists k\in\Z,\; \sum_{n\in\supp(\tau(h))} \chi_I(n-k) \equiv 1\pmod 2. \]
	Note that, if the sum is even for all $k$, then $\chi_I$ satisfies an invertible linear recurrence relation modulo $2$ (note that $\supp(\tau(h))\ne\emptyset$), hence $\chi_I$ is periodic by the pigeonhole principle. We conclude that
	\begin{itemize}[leftmargin=5mm]
		\item if $I$ is not periodic, then there exists $k$ such that the sum is odd. 
		\item if $I$ is $p$-periodic, it suffices to check for $k=0,1,\ldots,p-1$.
	\end{itemize}
	\item If $\tau(h)\notin\bigoplus_\Z C_2$, then $\tau(chc^{-1})=\tau(h)$ implies that $\tau(c)\in C_{C_2\wr\Z}(h)$ which is cyclic \cite[Lemma 2.5]{cox2018degree}. We deduce the existence of $r\in G_I$, $m,n\in\Z$ and $p,q\in\{0,1\}$ such that $g= r^mz^p$ and $c=r^nz^q$. Since $z$ is central, we can compute $chc^{-1}=h$. It follows that $h\not\sim hz$.
\end{itemize}

We can now decide the general Conjugacy Problem in $G_I$:
\begin{algorithm}[H]
	\caption{(Is $g\sim h$? Case when $I$ recursive)}\label{alg:final}
	\begin{algorithmic}[1]
		\Require{$g,h\in G_I$}
		\State Answer $\gets$ No
		\If{$\tau(g)\sim\tau(h)$} 
			\State Find $c\in G_I$ such that $\tau(cgc^{-1})=\tau(h)$
			\If{$cgc^{-1}=h$} \Comment{Theorem \ref{thm:WP}}
				\State Answer $\gets$ Yes
			\ElsIf{$h\sim hz$} \Comment{Sub-routine}
				\State Answer $\gets$ Yes
			\EndIf
		\EndIf
		\State\Return Answer
	\end{algorithmic}
\end{algorithm}
For lines 2,3, we can use \cite{Sale} which solves the Conjugacy Search Problem.
\end{proof}
\begin{rem}
	Observe that the sub-routine works for non-recursive $I$. We only need that $I$ is recursive when checking if $cgc^{-1}=h$. This shows that the Word Problem and the Conjugacy Problem are Turing equivalent.
\end{rem}
\begin{ques}
	Does $G(\Z^d,I)$ have decidable Conjugacy Problem as soon as $I$ is recursive? Equivalently, does Lemma \ref{lem:Turing_conjugate} fail in $\Z^d$ for $1<d<\infty$?
\end{ques}



\section{Miscellaneous}

\subsection{Residual finiteness} We characterise pairs $(H,I)$ such that $G(H,I)$ is residually finite. Some ideas towards this result were already given in \cite{ChapmanLubotzky}.

\begin{prop} \label{prop:residually_finite}
	$G(H,I)$ is residually finite if and only if $H$ is residually finite and $I$ is a union of cosets of a finite-index (normal) subgroup $K\le_{\mathrm{f.i.}}H$.
\end{prop}
\begin{proof}
	We first suppose that $G(H,I)$ is residually finite. Since $H\le G(H,I)$, we deduce that $H$ is residually finite too. Moreover, there exists $\varphi\colon G(H,I)\to F$ with $F$ finite and $\varphi(z)\ne 1_F$. Let $K=\ker\varphi|_H$. Since $F$ is finite, $K$ has finite-index inside $H$. Moreover, for all $k\in K$, we have
	\begin{align*}
		h\in I
		& \iff \varphi([a,hah^{-1}])\ne 1_F \\
		& \iff \varphi([a,(hk)a(hk)^{-1}])\ne 1_F \\
		& \iff hk\in I,
	\end{align*}
	proving that $I=IK$, i.e., $I$ is a union of cosets of $K$.
	
	\bigskip
	
	Reciprocally, if $H$ is residually finite, then $C_2\wr H$ is residually finite \cite{Grunberg}. It follows that, for each $g\ne 1,z$, there exists an epimorphism $\varphi\colon C_2\wr\Z\to F$ with $F$ finite such that $\varphi(\tau(g))\ne 1_F$. This reduces the question to $g=z$.
	
	Finally, if $I$ is a union of cosets of $K$ then there exists $J\subset G/K$ such that $I=\psi^{-1}(J)$ for $\psi\colon G\onto G/K$. Corollary \ref{lem:natural}(b) gives an epimorphism $G(H,I)\onto G(H/K,J)$ mapping $z\mapsto z\ne 1$ with $G(H/K,J)$ finite.
\end{proof}

\begin{ques}
	Is $G(H/K,J)$ the smallest quotient witnessing $z\ne 1$?
\end{ques}
\subsection{Context-free co-Word Problem} We have already seen that the Word Problem of $G_I$ can get much worse than that of $C_2\wr\Z$, when $I$ is not recursive. We study when the complexity does \emph{not} increase:

\begin{prop} Let $I\subset\Z$ be a symmetric set not containing $0$.
	\begin{enumerate}[leftmargin=8mm, label={\normalfont(\alph*)}]
		\item If $I$ is periodic, then $\coWP(G_I)$ is context-free (and $G_I$ embeds in $V$).
		\item If $\coWP(G_I)$ is context-free, then $I$ is eventually periodic.
	\end{enumerate}
\end{prop}
Here, we say that $I$ is \emph{eventually periodic} if there exists $p,K>0$ such that, for all $i>K$, we have $i\in I$ if and only if $i+p\in I$.
\begin{proof}
(a) $G_I$ is residually finite (Proposition \ref{prop:residually_finite}), hence there exists a finite-index subgroup $H\le G_I$ such that $z\notin H$. In particular, $\tau|_H$ provides an isomorphism between $H$ and a finite-index subgroup of $C_2\wr\Z$.

We now recall that $C_2\wr\Z$ embeds in $V$, and that the property \say{embeddable in $V$} passes to subgroups and finite-index overgroups. We deduce that $H$ and $G_I$ embed in Thomson's $V$, hence $\coWP(G_I)$ is context-free.

\medskip

(b) Since $\coWP(G_I)$ is context-free, the solution set of any instance of the Knapsack problem must be semi-linear \cite{Knapsack_coCF, Knapsack_semilinear}. In particular
\[ \{ (n_1,n_2,n_3,n_4)\in\N^4 \mid at^{n_1}at^{-n_2}at^{n_3}at^{-n_4}=1\} = \{(n,n,n,n)\mid n\in\N\setminus I\}\]
is semi-linear, meaning that $I$ is eventually periodic.
\end{proof}

\begin{conj}
	$\coWP(G_I)$ is context-free if and only if $I$ is periodic.
\end{conj}

\subsection{The Isomorphism Problem} It is known that $G_I\simeq G_J$ if and only if $I=J$ \cite{394830,fournier2023no}. We extend this result to the classes $G(H,I)$ when the group $H$ satisfies Kaplansky's unit conjecture. Notably, this covers all left-orderable groups (including $H=\Z$) and groups with the unique product property.

\begin{prop} \label{prop:iso_prob}
	Consider a group $H$ such that $\F_2[H]$ only has trivial units, and two symmetric subsets  $I,J\subset H$ satisfying $I,J\not\ni 1_H$. Then
	\[ G(H,I)\simeq G(H,J) \iff \bigl( \exists \psi\in\Aut(H),\; I=\psi(J)\bigr). \]
\end{prop}
\begin{proof}
	The statement is true if $H=\{1_H\}$. Suppose $H$ is infinite torsion-free. Consider an isomorphism $\varphi\colon G(H,I)\to G(H,J)$.
	\begin{itemize}[leftmargin=8mm]
		\item The center is $\Zc(G(H,I))=\{1,z\}$ hence $\varphi(z)=z$. It follows that $\varphi$ descends to an automorphism $\bar\varphi\colon C_2\wr H\to C_2\wr H$.
		\item The set of torsion elements in $G(H,I)$ is exactly $N(H,I)$. We deduce that $N(H,I)$ is a fully characteristic subgroup. Moreover
		\[ N(H,I)= \varphi^{-1}\bigl(N(H,J)\bigr)\]
		hence $\varphi$ descends to an automorphism $\bar{\bar\varphi}\colon H\to H$.
	\end{itemize}
	We identify $\bigoplus_HC_2\simeq\F_2[H]$ and observe that $\lla g\rra_{C_2\wr\Z}=(g)$ for all $g\in\F_2[H]$, where $(g)$ is the left ideal generated by $g$. In particular, $(g)=\F_2[H]$ if and only if $g$ is a unit of $\F_2[H]$. Since $\F_2[H]$ only has trivial units, we deduce that $\bar\varphi(a)=hah^{-1}$ for some $h\in H$. More generally, for all $k\in K$, we have
	\[ \bar\varphi\bigl(kak^{-1}\bigr)=\bigl(\bar{\bar\varphi}(k)h\bigr)a\bigl(\bar{\bar\varphi}(k)h\bigr)^{-1}. \]
	We have $k\in I$ if and only if $[a,kak^{-1}]=z$ in $G(H,I)$, if and only if $\varphi\bigl([a,kak^{-1}]\bigr)=z$ in $G(H,J)$. Let us compute this last expression:
	\begin{align*}
		\varphi\bigl([a,kak^{-1}]\bigr)
		& = \bigl[hah^{-1}z^\alpha,\bigl(\bar{\bar\varphi}(k)h\bigr)a\bigl(\bar{\bar\varphi}(k)h\bigr)^{-1}z^\beta\bigr] \\
		& = \bigl[a, \bigl(h^{-1}\bar{\bar\varphi}(k)h\bigr)a\bigl(h^{-1}\bar{\bar\varphi}(k)h\bigr)^{-1}\bigr] 
	\end{align*}
	(with $\alpha,\beta\in\{0,1\}$). We conclude that $k\in I\Leftrightarrow \psi(k)\coloneqq h^{-1}\bar{\bar\varphi}(k)h\in J$.
\end{proof}
\begin{rem}
	Another way to find a continuum of isomorphism types (for a fixed finitely generated, infinite group $H$) is to look at Turing complexity:
	\begin{adjustwidth}{3mm}{3mm}
		\textbf{Claim.} If $WP(H,S)\subset S^\star$ has Turing degree $\mathbf a$ and $I\subset S^\star$ has Turing degree $\mathbf b$, then $WP(G(H,I))$ has Turing degree $\mathbf a\oplus\mathbf b$.	
	\end{adjustwidth}
	Since there exists a continuum of Turing degrees $\mathbf b\ge_T\mathbf a$, there exists a continuum of Turing degrees among $WP(G(H,I))$, hence a continuum of isomorphism types among the central extensions $G(H,I)$.
\end{rem}

\appendix

\section{Conjugacy Geodesics in $C_2\wr F_r$.}

In this appendix, we prove that $\ConjGeo(C_2\wr F_r,T)$ is context-free. We first extract a description of elements $g=(\bfu,h)\in C_2\wr F_r$ which have minimal length in their conjugacy class from \cite[\S3]{mercier_arxiv}.

We first treat the case $h\ne 1$. We write $h=wcw^{-1}$ with $wcw^{-1}$ freely reduced and $c=c_1\ldots c_k$ cyclically reduced. We define
\[ \Sc_i=\{v\in F_r\mid \text{the reduced word for $v$ does not start with $c_{i-1}^{-1}$ or $c_i$}\}, \]
where $c_0\coloneqq c_k$ and $c_{k+1}\coloneqq c_1$. Finally, we define
\[\Tc_i = \Conv\Bigl(\!\bigl( (wc_1\ldots c_i)^{-1}\supp(\bfu) \cup \{1\}\bigr) \cap \Sc_i\Bigr). \]
(See Figure \ref{fig:def}.) Here $\Conv(X)$ is the smallest subtree containing $X$, with the convention $\Conv(\emptyset)=\emptyset$. We also denote $[x,y]=\Conv\bigl(\{x,y\}\bigr)$. 
\begin{center}
	\begin{tikzpicture}[scale=1.6]
	\definecolor{col0}{RGB}{216,27,96}
	\definecolor{col1}{RGB}{255,193,7}
	\definecolor{col1t}{RGB}{255,170,7}
	\definecolor{col2}{RGB}{30,136,229}
	
	\begin{scope}[every node/.style={circle, draw=black, inner sep=1.5pt}]
		\node[circle, fill=black, inner sep=1.5pt, label={below right:\footnotesize$w$}] (p0) at (0,0) {};
		\node[circle, fill=col1, inner sep=1.5pt] (p1) at (1.05,.1) {};
		\node[circle, fill=black, inner sep=1.5pt] (p2) at (2.1,-.1) {};
		\node[circle, fill=black, inner sep=1.5pt] (p3) at (3.1,.05) {};
		\node[circle, fill=col1, inner sep=1.5pt] (p4) at (4,-.05) {};
		\node[circle, fill=black, inner sep=1.5pt] (p5) at (5,0) {};
	\end{scope}
	
	\draw[-latex, thick, col0] (-1,-.05) to node[above] {$c_k$} (p0);
	\node[col0] at (-.5,-0.025) {$\times$};
	\draw[-latex, thick] (p0) to node[above] {$c_1$} (p1);
	\draw[-latex, thick] (p1) to node[above] {$c_2$} (p2);
	\draw[-latex, thick] (p2) -- (p3);
	\draw[-latex, thick] (p3) -- (p4);
	\draw[-latex, thick] (p4) to node[above] {$c_k$} (p5);
	\draw[-latex, thick, col0] (p5) to node[above] {$c_1$} (6.05,.1);
	\node[col0] at (5.525,0.05) {$\times$};
	
	\begin{pgfonlayer}{background}
		\clip (-1.5,-1.95) rectangle (6.5,2.1);
		\draw[col2, fill=col2!20] (-.5,-2) to[out=75, in=-75] (-.53,2.2) to (.53,2.2) to[out=-105, in=105] (.5,-2) -- cycle;
		\draw[col2, fill=col2!20] (.55,-2) to[out=75, in=-75] (.52,2.2) to (1.58,2.2) to[out=-105, in=105] (1.55,-2) -- cycle;
		\draw[col2, fill=col2!20] (1.6,-2) to[out=75, in=-75] (1.57,2.2) to (2.63,2.2) to[out=-105, in=105] (2.6,-2) -- cycle;
		\draw[col2, fill=col2!20] (4.5,-2) to[out=75, in=-75] (4.47,2.2) to (5.53,2.2) to[out=-105, in=105] (5.5,-2) -- cycle;
	\end{pgfonlayer}
	
	\begin{scope}[every node/.style={circle, draw=black}]
		\node[inner sep=1.1pt, fill=black] (n01) at (0.05,.6) {};
		\node[inner sep=1.3pt, fill=col1] (n02) at (0.2,1.2) {};
		\node[inner sep=1.1pt, fill=black] (n03) at (-0.05,1.1) {};
		\node[inner sep=1.3pt, fill=col1] (n04) at (-0.25,1.7) {};

		\node[inner sep=1.1pt, fill=black, label={[label distance=-1.5mm]left:\footnotesize$1$}] (e0) at (.2,-1.7) {};
		\node[inner sep=1.1pt, fill=black] (e1) at (.1,-1.2) {};
		\node[inner sep=1.3pt, fill=col1] (e2) at (-.05,-.6) {};
		
		\draw[ultra thick] (n04) -- (n03) -- (n01) -- (p0) -- (e2);
		\draw[ultra thick] (n02) -- (n01);
		\draw[col1, thick] (n04) -- (n03) -- (n01) -- (p0) -- (e2);
		\draw[col1, thick] (n02) -- (n01);
		
		\draw[-latex, thick] (e0) -- (e1);
		\draw[-latex, thick] (e1) -- (e2);
		
		\node[inner sep=1.1pt, fill=black] (n11) at (1.1,.6) {};
		\node[inner sep=1.1pt, fill=black] (n12) at (1.25,1.2) {};
		\node[inner sep=1.3pt, fill=col1] (n14) at (1.15,1.75) {};
		\node[inner sep=1.1pt, fill=black] (n15) at (1,-.6) {};
		\node[inner sep=1.3pt, fill=col1] (n16) at (0.9,-1.2) {};
		\draw[ultra thick] (n16) -- (n15) -- (p1) -- (n11) -- (n12) -- (n14);
		\draw[col1, thick] (n16) -- (n15) -- (p1) -- (n11) -- (n12) -- (n14);
		
		\node[inner sep=1.1pt, fill=black] (n51) at (5.05,.6) {};
		\node[inner sep=1.3pt, fill=col1] (n52) at (4.95,1.1) {};
		\node[inner sep=1.1pt, fill=black] (n53) at (5.2,1.2) {};
		\node[inner sep=1.3pt, fill=col1] (n54) at (5.1,1.7) {};
		\node[inner sep=1.3pt, fill=col1] (n55) at (5.3,1.65) {};
		\node[inner sep=1.3pt, fill=col1] (n56) at (4.85,-1.2) {};
		
		\node[inner sep=1.1pt, fill=black, label={[label distance=-1.5mm]left:\footnotesize$h$}] (h0) at (5.2,-1.7) {};
		\node[inner sep=1.1pt, fill=black] (h1) at (5.1,-1.2) {};
		\node[inner sep=1.3pt, fill=col1] (h2) at (4.95,-.6) {};
		
		\draw[ultra thick] (n52) -- (n51) -- (p5) -- (h2) -- (n56);
		\draw[ultra thick] (n51) -- (n53) -- (n54);
		\draw[ultra thick] (n53) -- (n55);
		\draw[col1, thick] (n52) -- (n51) -- (p5) -- (h2) -- (n56);
		\draw[col1, thick] (n51) -- (n53) -- (n54);
		\draw[col1, thick] (n53) -- (n55);
		
		\draw[-latex, thick] (h2) -- (h1);
		\draw[-latex, thick] (h1) -- (h0);
	\end{scope}
	
	\node[col2] at (0,1.93) {\small$w\Sc_0$};
	\node[col2] at (1.05,1.93) {\small$wc_1\Sc_1$};
	\node[col2] at (2.1,1.93) {\small$wc_1c_2\Sc_2$};
	\node[col2] at (5,1.93) {\small$wc\Sc_k$};
	
	\node[col1t, inner sep=1pt] (lT0) at (-.7,1.3) {\small$w\Tc_0$};
	\draw (lT0) -- (-.2,1.4);
\end{tikzpicture}
\captionsetup{margin=12mm, font=small}

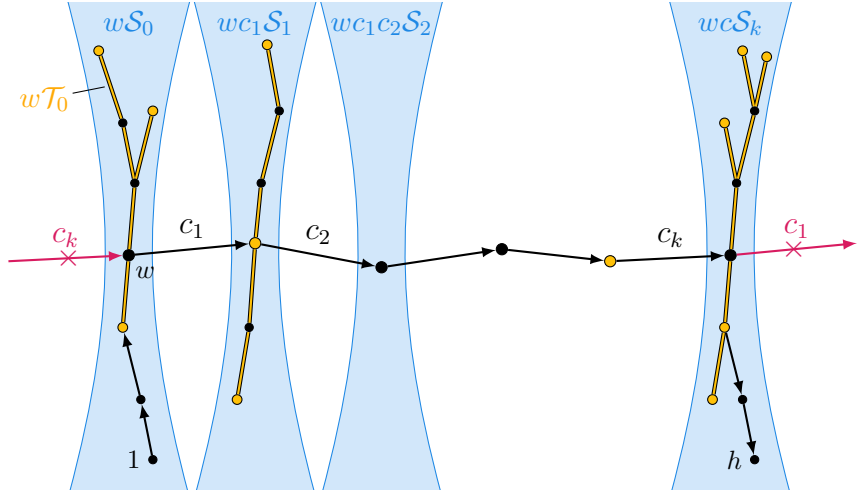
\captionof{figure}{An element $g=(\bfu,h)$ with $h\ne 1$, and the different associated constructions. Vertices in $\supp(\bfu)$ are marked in yellow.}
\label{fig:def}
\end{center}
\begin{prop} \label{prop:conj_min}
	Consider the generating set $T=\{(\delta_1,1)\}\cup B^\pm$, where $B$ is a basis of $F_r$. An element $g=(\bfu, h)\in C_2\wr F_r$ has minimal length in its conjugacy class if and only if it satisfies the following conditions:
	\begin{enumerate}[leftmargin=8mm, label={\normalfont(\arabic*)}]
		\item $\supp(\bfu)\subseteq w\Sc_0\sqcup wc_1\Sc_1 \sqcup wc_1c_2\Sc_2 \sqcup \ldots \sqcup wc_1\ldots c_k\Sc_k$,
		\item $\Tc_0\cap \Tc_k\subseteq [1,w^{-1}]$,
		\item $w^{-1}\supp(\bfu)\cap (wc)^{-1}\supp(\bfu) \cap [1,w^{-1}]=\emptyset$, and
		\item $w=1$ or $w^{-1}\in \Tc_0\cup \Tc_k$.
	\end{enumerate}	
\end{prop}
\begin{proof}
	We first explain the \say{only if} part. The key observation is
	\[ (\delta_x,1)(\bfu,h)(\delta_x,1)^{-1} = (\bfu+\delta_x-\delta_{hx},h)\]
	hence we can \say{transfer} lamps along $\la h\ra$ cosets. In particular,
	\begin{enumerate}[leftmargin=8mm, label=(\arabic*)]
		\item Since $w\Sc_0\sqcup wc_1\Sc_1 \sqcup wc_1c_2\Sc_2 \sqcup \ldots \sqcup wc_1\ldots c_{k-1}\Sc_{k-1}$ is a fundamental domain of the action $\la h\ra\acts F_r$, any element $(\bfu,h)$ is conjugated to an element $(\bfv,h)$ with
		\[ \supp(\bfv)\subseteq w\Sc_0\sqcup wc_1\Sc_1 \sqcup wc_1c_2\Sc_2 \sqcup \ldots \sqcup wc_1\ldots c_{k-1}\Sc_{k-1}. \]
		Moreover, if $g$ doesn't satisfy (1), this element will be strictly shorter (since we avoid some travelling along the axis $w\!\la c\ra$).
		
		\item If $g$ doesn't satisfy (2), we can move the lamps spanning $wc\Tc_k$ into $w\Sc_0$. The conjugated element  $g'$ is strictly shorter, since we only need to do one back-and-forth along the red edges on Figure \ref{fig:g_fails_2}, instead of two.
		\begin{center}
			\begin{tikzpicture}[scale=1.65]
	\definecolor{col0}{RGB}{216,27,96}
	\definecolor{col1}{RGB}{255,193,7}
	\definecolor{col1t}{RGB}{255,170,7}
	\definecolor{col2}{RGB}{30,136,229}
	
	\begin{scope}[every node/.style={circle, draw=black, inner sep=1.5pt}]
		\node[circle, fill=black, inner sep=1.5pt, label={below right:\footnotesize$w$}] (p0) at (0,0) {};
		\node[circle, fill=col1, inner sep=1.5pt] (p1) at (1.05,.1) {};
		\node[circle, fill=black, inner sep=1.5pt] (p2) at (2.1,-.1) {};
		\node[circle, fill=black, inner sep=1.5pt] (p3) at (3.1,.05) {};
		\node[circle, fill=col1, inner sep=1.5pt] (p4) at (4,-.05) {};
		\node[circle, fill=black, inner sep=1.5pt] (p5) at (5,0) {};
	\end{scope}
	
	\draw[-latex, thick] (p0) to node[above] {$c_1$} (p1);
	\draw[-latex, thick] (p1) to node[above] {$c_2$} (p2);
	\draw[-latex, thick] (p2) -- (p3);
	\draw[-latex, thick] (p3) -- (p4);
	\draw[-latex, thick] (p4) to node[above] {$c_k$} (p5);
	
	\begin{pgfonlayer}{background}
		\clip (-1.5,-1.95) rectangle (6.5,2.1);
		\draw[col2, fill=col2!20] (-.5,-2) to[out=75, in=-75] (-.53,2.2) to (.53,2.2) to[out=-105, in=105] (.5,-2) -- cycle;
		\draw[col2, fill=col2!20] (.55,-2) to[out=75, in=-75] (.52,2.2) to (1.58,2.2) to[out=-105, in=105] (1.55,-2) -- cycle;
		\draw[col2, fill=col2!20] (1.6,-2) to[out=75, in=-75] (1.57,2.2) to (2.63,2.2) to[out=-105, in=105] (2.6,-2) -- cycle;
		\draw[col2, fill=col2!20] (4.5,-2) to[out=75, in=-75] (4.47,2.2) to (5.53,2.2) to[out=-105, in=105] (5.5,-2) -- cycle;
	\end{pgfonlayer}
	
	\begin{scope}[every node/.style={circle, draw=black}]
		\node[inner sep=1.1pt, fill=black] (n01) at (0.05,.6) {};
		\node[inner sep=1.3pt, fill=col1] (n02) at (0.2,1.2) {};
		\node[inner sep=1.3pt, fill=col1] (n03) at (-0.05,1.1) {};
		\node[inner sep=1.3pt, fill=col1] (n04) at (-0.25,1.7) {};
		\node[inner sep=1.3pt, fill=col1] (n05) at (0.1,1.7) {};
		\node[inner sep=1.3pt, fill=col1] (n06) at (0.3,1.65) {};
		
		\node[inner sep=1.1pt, fill=black, label={[label distance=-1.5mm]left:\footnotesize$1$}] (e0) at (.2,-1.7) {};
		\node[inner sep=1.1pt, fill=black] (e1) at (.1,-1.2) {};
		\node[inner sep=1.3pt, fill=col1] (e2) at (-.05,-.6) {};
		
		\draw[ultra thick] (n04) -- (n03) -- (n01) -- (p0) -- (e2);
		\draw[ultra thick] (n01) -- (n02);
		\draw[ultra thick] (n05) -- (n02) -- (n06);
		\draw[col1, thick] 	(n04) -- (n03)
		(n05) -- (n02) -- (n06)
		(p0) -- (e2);
		\draw[col0, thick] 	(n01) -- (n02)
		(n03) -- (n01) -- (p0);
		
		\draw[-latex, thick] (e0) -- (e1);
		\draw[-latex, thick] (e1) -- (e2);
		
		\node[inner sep=1.1pt, fill=black] (n11) at (1.1,.6) {};
		\node[inner sep=1.1pt, fill=black] (n12) at (1.25,1.2) {};
		\node[inner sep=1.3pt, fill=col1] (n14) at (1.15,1.75) {};
		\node[inner sep=1.1pt, fill=black] (n15) at (1,-.6) {};
		\node[inner sep=1.3pt, fill=col1] (n16) at (0.9,-1.2) {};
		\draw[ultra thick] (n16) -- (n15) -- (p1) -- (n11) -- (n12) -- (n14);
		\draw[col1, thick] (n16) -- (n15) -- (p1) -- (n11) -- (n12) -- (n14);
		
		\node[inner sep=1.1pt, fill=black] (n51) at (5.05,.6) {};
		\node[inner sep=1.1pt, fill=black] (n52) at (4.95,1.1) {};
		\node[inner sep=1.1pt, fill=black] (n53) at (5.2,1.2) {};
		\node[inner sep=1.1pt, fill=black] (n54) at (5.1,1.7) {};
		\node[inner sep=1.1pt, fill=black] (n55) at (5.3,1.65) {};
		\node[inner sep=1.3pt, fill=col1] (n56) at (4.85,-1.2) {};
		
		\node[inner sep=1.1pt, fill=black, label={[label distance=-1.5mm]left:\footnotesize$h$}] (h0) at (5.2,-1.7) {};
		\node[inner sep=1.1pt, fill=black] (h1) at (5.1,-1.2) {};
		\node[inner sep=1.3pt, fill=col1] (h2) at (4.95,-.6) {};
		
		\draw[ultra thick] (p5) -- (h2) -- (n56);
		\draw[col1, thick] (p5) -- (h2) -- (n56);
		\draw[dotted,thick] (p5) -- (n51) -- (n52);
		\draw[dotted, thick] (n51) -- (n53) -- (n54);
		\draw[dotted, thick] (n53) -- (n55);
		
		\draw[-latex, thick] (h2) -- (h1);
		\draw[-latex, thick] (h1) -- (h0);
	\end{scope}
	
	\node[col2] at (0,1.93) {\small$w\Sc_0$};
	\node[col2] at (1.05,1.93) {\small$wc_1\Sc_1$};
	\node[col2] at (2.1,1.93) {\small$wc_1c_2\Sc_2$};
	\node[col2] at (5,1.93) {\small$wc\Sc_k$};
	
\end{tikzpicture}
\captionsetup{margin=12mm, font=small}
\captionof{figure}{A shorter element $g'\sim g$, using the failure of (2).}
\label{fig:g_fails_2}
		\end{center}
		
		\item If $g$ doesn't satisfy (3), we can cancel out some lamps in $\supp(\bfu)\cap [1,w]$ and $\supp(\bfu)\cap [h,hw]$ in the same $\la h\ra$ orbit, marked in red on Figure \ref{fig:g_fails_3}.
		\begin{center}
			\begin{tikzpicture}[scale=1.65]
	\definecolor{col0}{RGB}{216,27,96}
	\definecolor{col1}{RGB}{255,193,7}
	\definecolor{col1t}{RGB}{255,170,7}
	\definecolor{col2}{RGB}{30,136,229}
	
	\begin{scope}[every node/.style={circle, draw=black, inner sep=1.5pt}]
		\node[circle, fill=black, inner sep=1.5pt, label={below right:\footnotesize$w$}] (p0) at (0,0) {};
		\node[circle, fill=col1, inner sep=1.5pt] (p1) at (1.05,.1) {};
		\node[circle, fill=black, inner sep=1.5pt] (p2) at (2.1,-.1) {};
		\node[circle, fill=black, inner sep=1.5pt] (p3) at (3.1,.05) {};
		\node[circle, fill=col1, inner sep=1.5pt] (p4) at (4,-.05) {};
		\node[circle, fill=black, inner sep=1.5pt] (p5) at (5,0) {};
	\end{scope}
	
	\draw[-latex, thick] (p0) to node[above] {$c_1$} (p1);
	\draw[-latex, thick] (p1) to node[above] {$c_2$} (p2);
	\draw[-latex, thick] (p2) -- (p3);
	\draw[-latex, thick] (p3) -- (p4);
	\draw[-latex, thick] (p4) to node[above] {$c_k$} (p5);
	
	\begin{pgfonlayer}{background}
		\clip (-1.5,-1.95) rectangle (6.5,2.1);
		\draw[col2, fill=col2!20] (-.5,-2) to[out=75, in=-75] (-.53,2.2) to (.53,2.2) to[out=-105, in=105] (.5,-2) -- cycle;
		\draw[col2, fill=col2!20] (.55,-2) to[out=75, in=-75] (.52,2.2) to (1.58,2.2) to[out=-105, in=105] (1.55,-2) -- cycle;
		\draw[col2, fill=col2!20] (1.6,-2) to[out=75, in=-75] (1.57,2.2) to (2.63,2.2) to[out=-105, in=105] (2.6,-2) -- cycle;
		\draw[col2, fill=col2!20] (4.5,-2) to[out=75, in=-75] (4.47,2.2) to (5.53,2.2) to[out=-105, in=105] (5.5,-2) -- cycle;
	\end{pgfonlayer}
	
	\begin{scope}[every node/.style={circle, draw=black}]
		\node[inner sep=1.1pt, fill=black] (n01) at (0.05,.6) {};
		\node[inner sep=1.3pt, fill=col1] (n02) at (0.2,1.2) {};
		\node[inner sep=1.3pt, fill=col1] (n03) at (-0.05,1.1) {};
		\node[inner sep=1.3pt, fill=col1] (n04) at (-0.25,1.7) {};
		\node[inner sep=1.3pt, fill=col1] (n05) at (0.1,1.7) {};
		\node[inner sep=1.3pt, fill=col1] (n06) at (0.3,1.65) {};
		
		\node[inner sep=1.1pt, fill=black, label={[label distance=-1.5mm]left:\footnotesize$1$}] (e0) at (.2,-1.7) {};
		\node[inner sep=1.1pt, fill=black] (e1) at (.1,-1.2) {};
		\node[inner sep=3pt, draw=none, fill=col0, opacity=.5] at (-0.05,-.6) {};
		\node[inner sep=1.1pt, fill=black] (e2) at (-0.05,-.6) {};
		
		\draw[ultra thick] (n04) -- (n03) -- (n01) -- (p0);
		\draw[ultra thick] (n01) -- (n02);
		\draw[ultra thick] (n05) -- (n02) -- (n06);
		\draw[col1, thick] 	(n04) -- (n03) -- (n01) -- (p0)
		(n02) -- (n01)
		(n05) -- (n02) -- (n06);
		
		\draw[-latex, thick] (e0) -- (e1);
		\draw[-latex, thick] (e1) -- (e2);
		\draw[-latex, thick] (e2) -- (p0);
		
		\node[inner sep=1.1pt, fill=black] (n11) at (1.1,.6) {};
		\node[inner sep=1.1pt, fill=black] (n12) at (1.25,1.2) {};
		\node[inner sep=1.3pt, fill=col1] (n14) at (1.15,1.75) {};
		\node[inner sep=1.1pt, fill=black] (n15) at (1,-.6) {};
		\node[inner sep=1.3pt, fill=col1] (n16) at (0.9,-1.2) {};
		\draw[ultra thick] (n16) -- (n15) -- (p1) -- (n11) -- (n12) -- (n14);
		\draw[col1, thick] (n16) -- (n15) -- (p1) -- (n11) -- (n12) -- (n14);
		
		\node[inner sep=1.3pt, fill=col1] (n56) at (4.85,-1.2) {};
		
		\node[inner sep=1.1pt, fill=black, label={[label distance=-1.5mm]left:\footnotesize$h$}] (h0) at (5.2,-1.7) {};
		\node[inner sep=1.1pt, fill=black] (h1) at (5.1,-1.2) {};
		\node[inner sep=3pt, draw=none, fill=col0, opacity=.5] at (4.95,-.6) {};
		\node[inner sep=1.1pt, fill=black] (h2) at (4.95,-.6) {};
		
		\draw[ultra thick] (p5) -- (h2) -- (n56);
		\draw[col1, thick] (p5) -- (h2) -- (n56);
		
		\draw[-latex, thick] (h2) -- (h1);
		\draw[-latex, thick] (h1) -- (h0);
	\end{scope}
	
	\node[col2] at (0,1.93) {\small$w\Sc_0$};
	\node[col2] at (1.05,1.93) {\small$wc_1\Sc_1$};
	\node[col2] at (2.1,1.93) {\small$wc_1c_2\Sc_2$};
	\node[col2] at (5,1.93) {\small$wc\Sc_k$};
	
\end{tikzpicture}
\captionsetup{margin=12mm, font=small}
\captionof{figure}{A shorter element $g''\sim g'$, using the failure of (3).}
\label{fig:g_fails_3}
		\end{center}
		
		\item Finally, if $g$ doesn't satisfy (4), we can consider another conjugate. Write $w$ as a freely reduced word $w_1w_2\ldots w_\ell=w_1\tilde w$. The conjugate element
		\[ (\mathbf 0, w_1)^{-1}(\bfu, h)(\mathbf 0,w_1)=(w_1^{-1}\cdot\bfu, \tilde wc\tilde w^{-1}) \]
		is strictly shorter than $g$, since we avoid the red segments on Figure \ref{fig:g_fails_4}.
		\begin{center}
			\begin{tikzpicture}[scale=1.6]
	\definecolor{col0}{RGB}{216,27,96}
	\definecolor{col1}{RGB}{255,193,7}
	\definecolor{col1t}{RGB}{255,170,7}
	\definecolor{col2}{RGB}{30,136,229}
	
	\begin{scope}[every node/.style={circle, draw=black, inner sep=1.5pt}]
		\node[circle, fill=black, inner sep=1.5pt, label={below right:\footnotesize$\tilde w$}] (p0) at (0,0) {};
		\node[circle, fill=col1, inner sep=1.5pt] (p1) at (1.05,.1) {};
		\node[circle, fill=black, inner sep=1.5pt] (p2) at (2.1,-.1) {};
		\node[circle, fill=black, inner sep=1.5pt] (p3) at (3.1,.05) {};
		\node[circle, fill=col1, inner sep=1.5pt] (p4) at (4,-.05) {};
		\node[circle, fill=black, inner sep=1.5pt] (p5) at (5,0) {};
	\end{scope}
	
	\draw[-latex, thick] (p0) to node[above] {$c_1$} (p1);
	\draw[-latex, thick] (p1) to node[above] {$c_2$} (p2);
	\draw[-latex, thick] (p2) -- (p3);
	\draw[-latex, thick] (p3) -- (p4);
	\draw[-latex, thick] (p4) to node[above] {$c_k$} (p5);
	
	\begin{pgfonlayer}{background}
		\clip (-1.5,-1.95) rectangle (6.5,2.1);
		\draw[col2, fill=col2!20] (-.5,-2) to[out=75, in=-75] (-.53,2.2) to (.53,2.2) to[out=-105, in=105] (.5,-2) -- cycle;
		\draw[col2, fill=col2!20] (.55,-2) to[out=75, in=-75] (.52,2.2) to (1.58,2.2) to[out=-105, in=105] (1.55,-2) -- cycle;
		\draw[col2, fill=col2!20] (1.6,-2) to[out=75, in=-75] (1.57,2.2) to (2.63,2.2) to[out=-105, in=105] (2.6,-2) -- cycle;
		\draw[col2, fill=col2!20] (4.5,-2) to[out=75, in=-75] (4.47,2.2) to (5.53,2.2) to[out=-105, in=105] (5.5,-2) -- cycle;
	\end{pgfonlayer}
	
	\begin{scope}[every node/.style={circle, draw=black}]
		\node[inner sep=1.1pt, fill=black] (n01) at (0.05,.6) {};
		\node[inner sep=1.3pt, fill=col1] (n02) at (0.2,1.2) {};
		\node[inner sep=1.3pt, fill=col1] (n03) at (-0.05,1.1) {};
		\node[inner sep=1.3pt, fill=col1] (n04) at (-0.25,1.7) {};
		\node[inner sep=1.3pt, fill=col1] (n05) at (0.1,1.7) {};
		\node[inner sep=1.3pt, fill=col1] (n06) at (0.3,1.65) {};
		
		\node[inner sep=1.1pt, fill=black] (e0) at (.2,-1.7) {};
		\node[inner sep=1.1pt, fill=black, label={[label distance=-1.5mm]left:\footnotesize$1$}] (e1) at (.1,-1.2) {};
		\node[inner sep=1.1pt, fill=black] (e2) at (-0.05,-.6) {};
		
		\draw[ultra thick] (n04) -- (n03) -- (n01) -- (p0);
		\draw[ultra thick] (n01) -- (n02);
		\draw[ultra thick] (n05) -- (n02) -- (n06);
		\draw[col1, thick] 	(n04) -- (n03) -- (n01) -- (p0)
		(n02) -- (n01)
		(n05) -- (n02) -- (n06);
		
		\draw[col0, ultra thick, opacity=.5] (e0) -- (e1);
		\draw[-latex, thick] (e1) -- (e2);
		\draw[-latex, thick] (e2) -- (p0);
		
		\node[inner sep=1.1pt, fill=black] (n11) at (1.1,.6) {};
		\node[inner sep=1.1pt, fill=black] (n12) at (1.25,1.2) {};
		\node[inner sep=1.3pt, fill=col1] (n14) at (1.15,1.75) {};
		\node[inner sep=1.1pt, fill=black] (n15) at (1,-.6) {};
		\node[inner sep=1.3pt, fill=col1] (n16) at (0.9,-1.2) {};
		\draw[ultra thick] (n16) -- (n15) -- (p1) -- (n11) -- (n12) -- (n14);
		\draw[col1, thick] (n16) -- (n15) -- (p1) -- (n11) -- (n12) -- (n14);
		
		\node[inner sep=1.3pt, fill=col1] (n56) at (4.85,-1.2) {};
		
		\node[inner sep=1.1pt, fill=black] (h0) at (5.2,-1.7) {};
		\node[inner sep=1.1pt, fill=black, label={[label distance=-1mm]right:\footnotesize$\tilde wc\tilde w^{-1}$}] (h1) at (5.1,-1.2) {};
		\node[inner sep=1.1pt, fill=black] (h2) at (4.95,-.6) {};
		
		\draw[ultra thick] (p5) -- (h2) -- (n56);
		\draw[col1, thick] (p5) -- (h2) -- (n56);
		
		\draw[-latex, thick] (h2) -- (h1);
		\draw[col0, ultra thick, opacity=.5] (h1) -- (h0);
	\end{scope}
	
	\node[col2] at (0,1.93) {\small$\tilde w\Sc_0$};
	\node[col2] at (1.05,1.93) {\small$\tilde wc_1\Sc_1$};
	\node[col2] at (2.1,1.93) {\small$\tilde wc_1c_2\Sc_2$};
	\node[col2] at (5,1.93) {\small$\tilde wc\Sc_k$};
\end{tikzpicture}
\captionsetup{margin=12mm, font=small}
\captionof{figure}{A shorter element $g'''\sim g''$, using the failure of (4).}
\label{fig:g_fails_4}
		\end{center}
	\end{enumerate} 
	
	For the \say{if} part, we observe that all the reductions made in \cite[\S3]{mercier_arxiv} to go from $g$ satisfying (1-4) to a conjugacy geodesic actually preserve the length.
\end{proof}
\begin{rem}
	The statement can be simplified in rank $r=1$ since $w=1$ and $\Sc_i=\{1\}$. We only require (1') $\supp(\bfu)\subset [1,h]$ and (3') $\supp(\bfu)\not\supset\{1,h\}$.
\end{rem}
\bigskip

Next we treat the easier case $h=1$, i.e., $g\in \bigoplus_{F_r}C_2$.
\begin{prop}
	An element $g=(\bfu,1)\in C_2\wr F_r$ has minimal length its conjugacy class if and only if $\bfu=\mathbf 0$ or $1\in\Conv(\supp(\bfu))$.
\end{prop}
\begin{proof}
	Let $E_\bfu$ be the number of edges in the subtree $\Conv(\supp(\bfu)\cup\{1\})$. We have $\norm{(\bfu,1)}_T=2E_\bfu+\abs{\supp(\bfu)}$. Moreover, different elements obtained by conjugating only differ by a translation of the support:
	\[ (\bfv,h)(\bfu,1)(\bfv,h)^{-1}=(\bfu',1) \quad\text{where }\supp(\bfu')=h\cdot\supp(\bfu).\]
	The conclusion should be clear from here.
\end{proof}
\begin{thm} \label{thm:conjgeo_wreath}
	$\ConjGeo(C_2\wr F_r,T)$ is unambiguously context-free.
\end{thm}
\begin{proof}
	We denote $a=(\delta_1,1)$. 
	
	\textbullet\ We first consider elements $(\bfu,1)$ with $\supp\bfu\ne\emptyset$ and
	\[ \supp(\bfu) \subset \left\{ x\in F_r \;\big|\; \text{the reduced word for }x \text{ starts by } s\right\}. \]
	We call such elements \say{excursions in the direction $s$}.  We define production rules over the alphabet of terminals $T$ and variables $\{E_s \mid s\in B^\pm\}$:
	\[ E_s \ot sX_1X_2\ldots X_\ell s^{-1} \]
	where $s\in B^\pm$, the $X_i$'s are distinct elements of $\{a\}\cup\{E_v \mid v\ne s^{-1}\}$ and $\ell\ge 1$. Observe that the language defined by these rules with $E_s$ as the starting symbol is the language of geodesics for excursions in direction $s$.
	
	\textbullet\ We consider elements $(\bfu,h)$ with $h$ starting by $s$ and ending by $t$ and
	\[ \supp(\bfu) \subset\left\{ x\in F_r \;\big|\; \text{the common prefix of }h \text{ and } x \text{ is not }\varepsilon\text{ nor }h\right\}\]
	(see Figure \ref{fig:bridge}). (In particular $h\ne 1$, but we can have $h=s$ if $s=t$.) We call such elements \say{bridge in the direction $(s,t)$}.  We define production rules over the alphabet of terminals $T$ and variables $\{F_{s,t} \mid s,t\in B^\pm\}$:
	\[ F_{s,t} \ot s X_1X_2\ldots X_\ell F_{u,t} \]
	where $s,t,u\in B^\pm$ with $u\ne s^{-1}$, the $X_i$'s are distinct elements of $\{a\}\cup \{E_v\mid v\ne s^{-1},u\}$ and $\ell\ge 0$, as well as the extra rule $F_{s,s}\ot s$. Observe that the language defined by these rules with $F_{s,t}$ as the starting symbol is the language of geodesics for bridges in direction $(s,t)$.
	\begin{center}
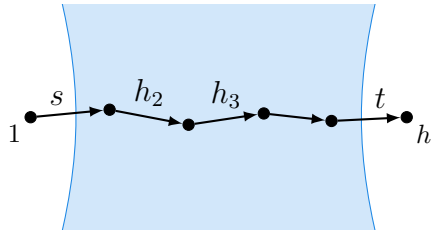

		\begin{tikzpicture}[scale=1]
	\definecolor{col0}{RGB}{216,27,96}
	\definecolor{col1}{RGB}{255,193,7}
	\definecolor{col1t}{RGB}{255,170,7}
	\definecolor{col2}{RGB}{30,136,229}
	
	\begin{scope}[every node/.style={circle, draw=black, inner sep=1.5pt}]
		\node[circle, fill=black, inner sep=1.5pt, label={below left:\footnotesize$1$}] (p0) at (0,0) {};
		\node[circle, fill=black, inner sep=1.5pt] (p1) at (1.05,.1) {};
		\node[circle, fill=black, inner sep=1.5pt] (p2) at (2.1,-.1) {};
		\node[circle, fill=black, inner sep=1.5pt] (p3) at (3.1,.05) {};
		\node[circle, fill=black, inner sep=1.5pt] (p4) at (4,-.05) {};
		\node[circle, fill=black, inner sep=1.5pt, label={below right:\footnotesize$h$}] (p5) at (5,0) {};
	\end{scope}
	
	\draw[-latex, thick] (p0) to (p1);
	\node at (.35,.25) {$s$};
	\draw[-latex, thick] (p1) to node[above] {$h_2$} (p2);
	\draw[-latex, thick] (p2) to node[above] {$h_3$} (p3);
	\draw[-latex, thick] (p3) -- (p4);
	\draw[-latex, thick] (p4) to (p5);
	\node at (4.65,.25) {$t$};
	
	\begin{pgfonlayer}{background}
		\clip (-1.5,-1.5) rectangle (6.5,1.5);
		\draw[col2, fill=col2!20] (.3,-2) to[out=75, in=-75] (.3,2) to (4.7,2) to[out=-105, in=105] (4.7,-2) -- cycle;
	\end{pgfonlayer}
\end{tikzpicture}
\captionsetup{margin=12mm, font=footnotesize}
\captionof{figure}{Set where $\bfu$ can possibly be supported.}
\label{fig:bridge}
	\end{center}

	\textbullet\ Conjugacy geodesics for elements of the form $(\bfu, h)$ with $h\ne 1$ cyclically reduced are produced by the rules
	\[ S \ot X_1\ldots X_k F_{s,t} Y_1\ldots Y_\ell \]
	where $s,t\in B^\pm$ with $s\ne t^{-1}$, the letters $X_i,Y_j$ are distinct elements of $\{a\}\cup\{ E_v\mid v\ne s,t^{-1}\}$ and $k,\ell\ge 0$.
	
	\textbullet\ Conjugacy geodesics for elements of the form $(\bfu, h)$ where $h\ne 1$ is not cyclically reducible are produced by the following rules:
	\begin{itemize}[leftmargin=8mm, label=-]
		\item $S \ot X_1\ldots X_k \cdot tS_t t^{-1}\cdot  Y_1\ldots Y_\ell$ where $t\in B^\pm$, the $X_i,Y_j$ are distinct elements of $\{a\}\cup\{E_v\mid v\ne t\}$ and $k+\ell\ge 1$.
		\item $S_s \ot X_1\ldots X_k \cdot tS_t t^{-1}\cdot  Y_1\ldots Y_\ell$ where $s,t\in B^\pm$ with $s\ne t^{-1}$, the  $X_i,Y_j$ are distinct elements of $\{a\}\cup\{E_v\mid v\ne s^{-1},t\}$ and $k,\ell\ge 0$.
		\item $S_s \ot X_1\ldots X_k \cdot F_{t,u} \cdot  Y_1\ldots Y_\ell $ where $s,t,u\in B^\pm$ with $s\ne t^{-1}\ne u$, the $X_i,Y_j$ are distinct elements of $\{a\}\cup\{E_v\mid v\ne s^{-1},t,u^{-1}\}$ and $k,\ell\ge 0$.
	\end{itemize}

	\textbullet\ Finally, conjugacy geodesics for elements of the form $(\bfu,1)$ are produced by the rules $S\ot \varepsilon$ (for $g=(\mathbf 0,1)$) and 
	\[ S \ot X_1X_2\ldots X_\ell \]
	where the $X_i$ are distinct elements of $\{a\}\cup\{E_v\mid v\in B^\pm\}$ and $\ell\ge 2$.
	
	We claim that the language produced when putting all these rules together, and considering $S$ as the starting symbol, is $\ConjGeo(C_2\wr F_2,T)$. Moreover, any conjugacy geodesic admits a unique leftmost derivation.
\end{proof}

\AtNextBibliography{\small}
\printbibliography

\end{document}